\documentclass[12pt]{amsart}
\usepackage{graphicx}

\newcommand{\p}{\partial}

\theoremstyle{plain}
      \newtheorem{thm}{Theorem}
      \newtheorem{cor}[thm]{Corollary}
      
      \newtheorem{lem}[thm]{Lemma}

      \newtheorem{prop}[thm]{Proposition}

\newcommand{\C}{\mathbb C}

\newcommand{\R}{\mathbb R}
\newcommand{\Z}{\mathbb Z}

\newcommand{\SLC}{SL(2,\mathbb C)}
\newcommand{\tL}{\tilde L}

\newcommand{\tX}{\tilde X}
\newcommand{\CP}{\mathbb{CP}}

\begin{document}

\title{Exceptional surgery and boundary slopes}
\author{Masaharu Ishikawa, Thomas W.\ Mattman and Koya Shimokawa}
\address{Department of Mathematics, Tokyo Metropolitan University,
                        1-1 Minami-Ohsawa, Hachioji-shi, Tokyo 192-0397,
Japan}
\email{mishika@math.matro-u.ac.jp}
\address{Department of Mathematics and Statistics,
                        California State University, Chico,
                        Chico CA 95929-0525, USA}
\email{TMattman@CSUChico.edu}

\address{Department of Mathematics,
                        Saitama University,
                    255 Shimo-Okubo Saitama-shi, 338-8570, Japan}
\email{kshimoka@rimath.saitama-u.ac.jp}
\thanks{This research is partially supported
by the Japan Society  for the Promotion of Science,
Grant-in-Aid for Scientific Research (B) 13440051.
The first author is partially supported by the
Japan Society  for the Promotion of Science for Young Scientists.
The third author is  partially supported
by the Ministry of Education, Culture, Sports, Science and Technology,
Grant-in-Aid for Young Scientists (B) 13740031}

\begin{abstract}
Let $X$ be a norm curve in the $\SLC$-character variety of a knot exterior
$M$. Let $t = \| \beta \| / \| \alpha \|$ be the ratio of the
Culler-Shalen  norms of two distinct non-zero classes $\alpha, \beta \in
H_1(\partial M,\Z)$. We demonstrate that either $X$ has exactly two associated strict
boundary slopes $\pm t$, or else there are strict boundary slopes $r_1$
and $r_2$ with $|r_1| > t$ and $|r_2| < t$. As a consequence, we show
that there are strict boundary slopes near cyclic, finite, and Seifert
slopes. We also prove that the diameter of the set of strict boundary 
slopes can be bounded below using the Culler-Shalen norm of those slopes.
\end{abstract}

\maketitle
\setlength{\baselineskip}{15pt}

\section{Introduction}

For a knot in a closed (i.e., compact and without boundary),
connected, orientable $3$-manifold,
let $M$ denote the exterior of the knot.
We fix a basis $(\mu, \lambda)$ of $H_1(\partial M,\mathbb Z)$.
The slope of $\gamma\in H_1(\p M,\mathbb Z)$ with respect to
this basis will be denoted $r_\gamma$.
That is, if $\gamma=a\mu +b\lambda$,
$r_\gamma = a/b \in \mathbb Q\cup \{1/0\}$.
  Let $M(r)$ denote {\em Dehn surgery} on a knot along slope $r$.
That is, $M(r)$ is the manifold obtained by attaching a solid torus $V$
to $M$ by a homeomorphism
of $\p V\to \p M$ which sends a meridian curve of $V$ to
a simple closed curve in $\p M$ of the given slope $r$.
If $\pi_1(M(r))$ is cyclic (respectively, finite),
we call $r$ a {\em cyclic} (resp., {\em finite}) {\em slope}.
If $M(r)$ admits the structure of a Seifert fibred space,
we call $r$ a {\em Seifert slope}.
An {\em essential surface} $F$ in $M$ is an incompressible and
orientable surface properly embedded in $M$,
no component of which is $\p$-parallel and no $2$-sphere component of
which bounds a $B^3$.
If the set $\{\p F\}$ is not empty,
it consists of a collection of parallel, simple closed curves in $\p M$.
We call the slope of such a curve obtained from an essential
surface $F$ a {\em boundary slope} and
we say that it is a {\em strict boundary slope} if we can choose
$F$ so that it is not the fibre of any fibration over the circle.

This paper deals with the connection
between boundary slopes and
cyclic, finite, and Seifert slopes.
If $M$ is hyperbolic, these last three types of slopes are examples
of exceptional slopes, i.e., $M(r)$ is not hyperbolic.
We also show that the diameter of the set of
strict boundary slopes can be bounded below in terms of the norms of
such slopes.

The set of characters of representations $\rho : \pi_1(M) \to \SLC$
can be identified with the points of a complex affine algebraic set $X(M)$,
which is called the {\em character variety} \cite{cs}.
For $\gamma \in\pi_1(M)$ we define the regular function
$I_\gamma : X(M) \to \mathbb C$
by $I_{\gamma}(\chi_\rho) = \text{trace}(\rho(\gamma))$.
By the Hurewicz isomorphism,
a class $\gamma \in L = H_1(\p M, \mathbb Z)$ determines an element
of $\pi_1(\p M) \subset \pi_1(M)$ well defined up to conjugacy.
A {\em norm curve} $X_0$ is a one-dimensional irreducible
component of $X(M)$ on which no $I_\gamma$
($\gamma \in L \setminus \{0\}$) is constant.
In this paper we will assume that $X(M)$ contains a norm curve.
For example, it is known that this assumption holds if $M$ is hyperbolic.

The terminology reflects the fact that we may associate
to $X_0$ a norm $\| \cdot \|_0$ on $H_1(\p M, \mathbb R)$
called a {\em Culler-Shalen norm} in the following manner.
Let $\tX_0$ be the smooth projective model of $X_0$,
which is birationally equivalent to $X_0$.
The birational map is regular at all but a finite number of points of
$\tX_0$, which are called {\em ideal points} of $\tX_0$.
The function $f_\gamma = I^2_\gamma -4$
is regular on $X_0$, and so can be pulled back to $\tX_0$.
We will also denote the pull-back by $f_\gamma$.
For $\gamma \in L$, the Culler-Shalen norm $\| \gamma \|_0$ is the
degree of $f_\gamma: \tX_0 \to \CP^1$.
The norm is extended to $H_1(\p M, \mathbb R)$ by linearity.

Fix a norm curve $X_0$ in the character variety $X(M)$
and denote by $\mathcal I$ the set of ideal points on $\tX_0$.
Let $s_0$ denote the minimal norm of $\|\cdot\|_0$,
i.e., $s_0=\text{min}\{\|\gamma\|_0; \gamma \in L, \gamma\neq 0\}$.
If $M$ is hyperbolic,
let $X_i$ denote a component of $X(M)$ which contains
the character of a discrete, faithful representation.
Note that $X_i$ is a norm curve by \cite[Proposition 1.1.1]{cgls}.
Let $\|\cdot\|_i$ denote the norm of $X_i$.
We define the {\em canonical norm} $\|\cdot\|_M$ on $H_1(M,\mathbb Z)$
to be the sum
\[
    \|\cdot\|_M=\|\cdot\|_1+\|\cdot\|_2+\cdots+\|\cdot\|_k
\]
as in \cite{bz2}.
Let $s_M$ denote the minimal norm of $\|\cdot\|_M$.

Note that $L = H_1(M,\mathbb Z)$ is a lattice of $V = H_1(M,\mathbb R)$,
i.e., a $\mathbb Z$-submodule of $V$ which is finitely generated
and spans $V$ as a vector space over $\mathbb R$.
Let $\tilde L$ denote a sublattice of $L$.
For an element $\gamma \in L$, let $\tilde \gamma$ denote
a primitive element in $\tilde L$ such that
$\tilde \gamma = q \gamma$ 
in $L$ for some $q \in \mathbb N$.
Let $s_{\tilde \gamma}$ denote the slope of $\tilde \gamma $
with respect to a basis $(\alpha, \beta)$ of $\tilde L$.

Now we state our main theorem.

\begin{thm}\label{main}
Let $M$ be a knot exterior
and $\alpha$ and $\beta$ be distinct non-zero elements in 
$L=H_1(\p M,\mathbb Z)$ which span a sublattice $\tilde L$.
Suppose $X(M)$ contains a norm curve
with the Culler-Shalen norm $\|\cdot\|_0$.
Then one of the following holds.
\begin{enumerate}
\item There are two distinct strict boundary classes $\gamma$ and $\delta$
such that $|s_{\tilde \gamma}|<{\|\beta\|_0}/{\|\alpha\|_0}$ and
$|s_{\tilde \delta}|> {\|\beta\|_0}/{\|\alpha\|_0}$.
In case $\alpha =\mu$ and $\beta=a\mu +b\lambda$ with $b>0$, we have
    $|r_\gamma -r_\beta|<{\|\beta\|_0}/{b\|\alpha\|_0}$ and
$|r_\delta -r_\beta|> {\|\beta\|_0}/{b\|\alpha\|_0}$.
\item There are exactly two distinct strict boundary classes $\gamma$
and $\delta$
associated to $\mathcal I$.
Moreover they satisfy $-s_{\tilde \gamma}={\|\beta\|_0}/{\|\alpha\|_0}$ and
$s_{\tilde \delta}={\|\beta\|_0}/{\|\alpha\|_0}$.
In case $\alpha =\mu$ and $\beta=a\mu +b\lambda$ with $b>0$,
$r_\beta -r_\gamma={\|\beta\|_0}/{b\|\alpha\|_0}$ and
$r_\delta -r_\beta={\|\beta\|_0}/{b\|\alpha\|_0}$.
\end{enumerate}

If $M$ is hyperbolic, the same statement also holds for
the canonical norm $\|\cdot\|_M$.
\end{thm}

As a direct corollary we have the following.

\begin{cor}\label{general}
Let $M$ be a knot exterior
and $\alpha$ and $\beta$ be distinct non-zero elements in 
$L=H_1(\p M,\mathbb Z)$ which span a sublattice $\tilde L$.
Suppose $X(M)$ contains a norm curve
with the Culler-Shalen norm $\|\cdot\|_0$.
Suppose $\beta$ has a particular property and
we have an upper bound $c$ on the norm of such a class,
i.e., $\|\beta\|_0\leq c$.
Then there is a strict boundary class $\gamma$ with
$|s_{\tilde \gamma}|\leq c/\|\alpha\|_0$.
In case $\alpha =\mu$ and $\beta=a\mu +b\lambda$ with $b>0$, we have
$|r_\gamma -r_\beta|\leq c/b\|\alpha\|_0$.

If $M$ is hyperbolic, the same statement also holds for
the canonical norm $\|\cdot\|_M$.
\end{cor}

Theorem \ref{main} and Corollary \ref{general} can be applied to
the study of relations between boundary slopes and 
cyclic, finite, or Seifert slopes.
A {\em small Seifert manifold} is a $3$-manifold which admits
the structure of a Seifert fibred space whose base orbifold is $S^2$
with at most three cone points.
A small Seifert manifold is irreducible if and only if it is not
$S^1 \times S^2$, and Haken if and only if it has infinite first
homology.

\begin{cor}\label{app}
Let $M$ be a knot exterior
and $\alpha$ and $\beta$ be distinct non-zero elements in
$L=H_1(\p M,\mathbb Z)$ which span a sublattice $\tilde L$.
Suppose $X(M)$ contains a norm curve.
\begin{enumerate}
\item
Suppose $\beta$ has minimal norm.
Then there is a strict boundary class $\gamma$ with
$|s_{\tilde \gamma}|\leq 1$.

Suppose further that $M$ is hyperbolic and $H^1(M,\mathbb Z_2)=\mathbb Z_2$.
Suppose $\alpha$ and $\beta$ are cyclic classes and are not
strict boundary classes.
Then there are two distinct strict boundary classes $\gamma$ and $\delta$
such that $|s_{\tilde \gamma}|<1$ and $|s_{\tilde \delta}|> 1$.
In case $\alpha =\mu$ and $\beta=a\mu +\lambda$, we have
$|r_\gamma -r_\beta|<1$ and $|r_\delta -r_\beta|> 1$.
\item
Suppose $M$ is hyperbolic and
$\beta$ is a finite class.
Then there is a strict boundary class $\gamma$ with
$|s_{\tilde \gamma}|\leq 3$.
In case $\alpha =\mu$ and $\beta=a\mu +b\lambda$ with $b>0$, we
 have $|r_\gamma -r_\beta|\leq 3/b$.
\item
Suppose that there is a class $\delta$ in $L$ such that
{\rm Hom}$(\pi_1(M(\delta)), PSL(2,\mathbb C))$ contains only
diagonalisable representations.
Suppose
$M(\beta)$ is an irreducible non-Haken small Seifert manifold.
Then there is a strict boundary class $\gamma$ with
$|s_{\tilde \gamma}|\leq 1+{2A}/{s_0}$,
where $A$ is the number of characters $\chi_\rho \in X_0$
of non-abelian representations $\rho \in R_0$ with $\rho(\beta)=\pm I$.
In case $\alpha =\mu$ and $\beta=a\mu +b\lambda$ with $b>0$, we have
$|r_\gamma -r_\beta|\leq (1+{2A}/s_0)/b$.
\end{enumerate}
\end{cor}

Note that if $M$ is the exterior of a (hyperbolic) knot in $S^3$
then it satisfies the conditions of Corollary~\ref{app},
i.e., $H^1(M,\mathbb Z_2)=\mathbb Z_2$ and
Hom$(\pi_1(S^3), PSL(2,\mathbb C)))$ contains only
diagonalisable representations.
Corollary~\ref{app}(1) improves a result of Dunfield\cite{d1} who showed
that for a cyclic slope $r$
there is a boundary slope in $(r-1,r+1)$.

If $M$ is the exterior of a hyperbolic knot in  a homotopy $3$-sphere,
we can take a preferred meridian-longitude pair for
$(\mu, \lambda)$. Then $r_\beta = a/b$ is the usual slope.
In this case, by~\cite[Theorem 1.1]{bz1},
$b$ in Corollary~\ref{app}(2) is either $1$ or $2$ and,
for the fillings in Corollary~\ref{app}(2) and (3), $b=1$ is
conjectured. (See Conjecture A in problem 1.77 of \cite{kirby}.)

Next we consider a relationship between the diameter of the set of
strict boundary slopes and the norms of these slopes.
Let $\mathcal B$ be the set of strict boundary slopes
 associated to
$\mathcal I$ with respect to a basis $(\mu,\lambda)$ of $L$.
As in \cite{cs3}, if $\infty \notin \mathcal B$,
let diam\,$\mathcal B$ denote the {\em diameter} of $\mathcal B$,
which is defined to be the difference between the greatest and
least elements of $\mathcal B$.
   From Theorem \ref{main} we obtain the following corollary.

\begin{cor}\label{diam}
Let $M$ be a knot exterior with $\infty \notin \mathcal B$.
Suppose $X(M)$ contains a norm curve $X_0$ with the norm $\|\cdot\|_0$.
Let $\beta$ be a strict boundary class associated to an ideal point
of $\mathcal I$ with $r_\beta = a/b$.
Then {\rm diam}\,$\mathcal B > \|\beta\|_0/b\|\mu\|_0$.
\end{cor}

Note that if $M$ is hyperbolic, by Proposition 9.1 of \cite{bz2}
we have
$\|\beta\|_0/b\|\mu\|_0=\|\beta\|_M/b\|\mu\|_M$.

In \cite{cs3} Culler and Shalen showed that if $M$ is the
exterior of a non-trivial, non-cable knot in an orientable
$3$-manifold with cyclic fundamental group,
then  {\rm diam}\,$\mathcal B \geq 2$.

The structure of the paper is as follows. In the next section
we give a brief introduction to the character variety and the
Culler-Shalen norm and how they apply to the study
of exceptional surgeries.
We prove Theorem 1 and Corollaries 2, 3, and 4 in Section~\ref{sec1}.
Finally, in Section~\ref{secex}
we discuss examples: the $(-2,3,n)$ pretzel knots and the twist knots.

\vskip 0.5cm

\section{Character variety, Culler-Shalen norm,
and exceptional surgery\label{secCS}}

In this section we recall the definitions  of
character varieties and Culler-Shalen norms.
The main references are the first chapter of~\cite{cgls} and~\cite{cs}.
Applications to finite surgery are developed in~\cite{bz1,bz2}.

Let $R(M)$ denote the set of representations $\rho : \pi_1(M) \to \SLC$.
It is easy to show that $R(M)$ is a complex affine algebraic set.
 The {\em character} of an element $\rho\in R(M)$ is the function
$\chi_\rho: \pi_1(M)\to\C$ defined by the trace map
$\chi_\rho(\gamma)=\text{trace}(\rho(\gamma))$.
The set of characters of the representations in $R(M)$
is also a complex affine algebraic set\cite{cs}.
We call it the {\em character variety} of $\pi_1(M)$
and denote it by $X(M)$.

For $\gamma \in \pi_1(M)$ we define the
regular function $I_\gamma : X(M) \to \mathbb C$
by $I_{\gamma}(\chi_\rho) = \text{trace}(\rho(\gamma))$.
The Hurewicz isomorphism induces an isomorphism
$H_1(\p M,\mathbb Z)\simeq \pi_1(\p M)\subset \pi_1(M)$. So
we can identify $L = H_1(\p M, \mathbb Z)$ with a subgroup
of $\pi_1(M)$, well defined up to conjugacy.
Thus each element $\gamma \in L$ determines a regular function $I_\gamma$.
A {\em norm curve} $X_0$ is a one-dimensional irreducible
component of $X(M)$ such that, for each $\gamma \in L \setminus \{0\}$,
$I_\gamma$ is not constant on $X_0$.
By~\cite[Proposition 2]{cs2},
any irreducible component of $X(M)$ containing the character of
a discrete, faithful representation is a norm curve of $X(M)$.
If $M$ is hyperbolic,
$X(M)$ contains the character of a discrete,
faithful representation, namely the holonomy representation, so there
will be at least one norm curve in the character variety.

The {\em Culler-Shalen norm} is a norm, on the
real vector space $H_1(\p M,\R)$,
associated to a norm curve $X_0$ in the following manner.
Let $f_\gamma=I^2_\gamma-4$. Since this function is regular,
it can be pulled back to $\tX_0$,
where $\tX_0$ is the smooth projective completion of $X_0$.
We again denote the pull-back by $f_\gamma$.
For $\gamma \in L$, define $\| \gamma \|_0$ to be the
degree of $f_\gamma: \tX_0 \to \CP^1$.
It is shown in~\cite[Proposition 1.1.2]{cgls} that
there exists a norm $\|\gamma\|_0$ on $H_1(\p M,\R)$ satisfying
(i) $\|\gamma\|_0=\text{degree} f_\gamma$ when $\gamma \in L$, and
(ii) the unit ball is a finite-sided polygon whose vertices
are rational multiples of strict boundary classes in $L$.
We call this norm the Culler-Shalen norm.


Let $\beta$ be a finite class which is not a strict boundary
class. Following a classification of Milnor~\cite{mi},
Boyer and Zhang~\cite{bz1} say that $\beta$ falls into one of six types
C, D, T, O, I, or Q. The notation refers to the fact that
$\pi_1(M(r))$ is an extension of a Cyclic, Dihedral, Tetrahedral, etc.\
group.

By~\cite[Corollary 1.1.4]{cgls},
a cyclic or C-type class
which is not a boundary class
realizes the minimal norm on $L = H_1(\p M, \Z)$.
In general, for a finite slope $r = r_{\beta}$ which is not a boundary
slope, $\beta$
realizes the minimal norm on a sublattice $\tL$ of some index $q$.
This is Proposition~9.3 of~\cite{bz2} which we restate here:
\vspace{1.5cm}

\begin{thm}[Proposition 9.3 \cite{bz2}]\label{finite}
Let $M$ be a hyperbolic knot exterior.
Let $s_M=\text{min}\{ \|\gamma\|_M; \gamma \in L =  H_1(\p M, \Z),\gamma
\neq 0\}$.
Suppose that $\beta$ is a finite class and not
a strict boundary class. 
Then there is an integer $q \in \{1, \ldots, 5\}$
and an index $q$ sublattice $\tL$ of $L$
such that $\| \beta \|_M \leq \| \gamma \|_M$ for all $0 \neq \gamma\in\tL$.

Moreover,
\begin{enumerate}
\item if $\beta$ is C-type, then $\|\beta\|_M=s_M$, i.e., $q = 1$;
\item if $\beta$ is D-type or Q-type, then $\|\beta\|_M\leq 2s_M$
         and $q \leq 2$;
\item if $\beta$ is T-type, then $\|\beta\|_M\leq s_M+4$ and $q \leq 3$;
\item if $\beta$ is I-type, then $\|\beta\|_M\leq s_M+8$ and
         $q \in \{ 1,2,3,5 \}$; and
\item {\rm (a)}if $\beta$ is O-type and $H_1 (M,\mathbb Z)$ has no
         non-trivial even torsion, then $\|\beta\|_M\leq s_M+6$ and
         $q \in \{2,4\}$ and
         {\rm (b)}if $\beta$ is O-type and $H_1 (M,\mathbb Z)$ has
         non-trivial even torsion, then $\|\beta\|_M\leq s_M+12$ and
         $q \leq 3$.
\end{enumerate}
\end{thm}

Note that there may be more than one choice of $q$ for a given finite
slope. For example, the O-type surgery $22$ of the $(-2,3,9)$ pretzel
knot realizes the minimal norm on sublattices of index $q=2$, $3$, and
$4$.

Next we refer to a result of Boyer and Ben Abdelghani.

\begin{thm}[Theorem C \cite{bb}]\label{bb}
Let $M$ be a knot exterior.
Suppose that there is
 a class $\delta$ in $L$ such that
{\rm Hom}$(\pi_1(M(\delta)), PSL(2,\mathbb C))$ contains only
diagonalisable representations.
Suppose, for a non-boundary class $\beta$,
$M(\beta)$ is an irreducible non-Haken small Seifert manifold.
Then $\| \beta \|_0=s_0 +2A$,
where $A$ is the number of characters $\chi_\rho \in X_0$
of non-abelian representations $\rho \in R_0$ with $\rho(\beta)=\pm I$.
\end{thm}

\vskip 0.5cm

\section{Proofs\label{sec1}}

Fix a norm curve $X_0$ in the character variety $X$
with set of ideal points $\mathcal I$.
Let $\Pi_x(f_\alpha)$ denote the order of the pole of $f_\alpha$
at $x$.

We start by stating the main tool of our proof.

\begin{prop}\label{prop}
Let $M$ be a knot exterior
and $\alpha$ and $\beta$ be elements in $L=H_1(\p M,\mathbb Z)$.
Suppose $X(M)$ contains a norm curve $X_0$
with the Culler-Shalen norm $\|\cdot\|_0$.
Then either
\begin{enumerate}
\item there are two distinct ideal points $x$ and $y$ such that
$\Pi_x(f_\alpha)/\|\alpha\|_0 < \Pi_x(f_\beta)/\|\beta\|_0$ and
$\Pi_y(f_\alpha)/\|\alpha\|_0 > \Pi_y(f_\beta)/\|\beta\|_0$,
or
\item for any ideal point $z$, we have
$\Pi_z(f_\alpha)/\|\alpha\|_0 = \Pi_z(f_\beta)/\|\beta\|_0$.
\end{enumerate}

If $M$ is hyperbolic, the same statement also holds for
the canonical norm $\|\cdot\|_M$.
\end{prop}

\begin{proof}
   From the definition of the norm $\|\cdot \|_0$, we have,
$\|\alpha\|_0=\sum_{x \in \mathcal I} \Pi_x(f_\alpha)$.
Hence we have
$\sum_{x \in \mathcal I}\Pi_z(f_\alpha)/\|\alpha\|_0 = \sum_{x \in
\mathcal I}\Pi_z(f_\beta)/\|\beta\|_0$.
Hence if (2) does not hold, then (1) holds.

The same argument applies to $\|\cdot\|_M$.
\end{proof}

\begin{lem}\label{two}
Suppose $M$ has a norm curve $X_0$ with the norm $\|\cdot\|_0$.
Then there are two ideal points in $\mathcal I$
whose associated strict boundary classes are distinct.
\end{lem}

\begin{proof}
Suppose, for a contradiction, that there is at most one
strict boundary class associated to $\mathcal I$.
If there is no strict boundary class associated to $\mathcal I$,
we have $\|\gamma\|_0=0$ for any element $\gamma \in L$.
If each strict boundary class associated to $\mathcal I$ is
equal to $\gamma \in L$, then we have $\|\gamma\|_0=0$.
In both cases we have contradictions to the fact that
$\|\cdot\|_0$ is a norm.
\end{proof}

\noindent
{\it Proof of Theorem $\ref{main}$.}
Suppose (1) of Proposition \ref{prop} holds.
Then we have ideal points $x$ and $y$ satisfying the inequalities
described in the proposition.

First suppose $\Pi_x(f_\alpha)>0$.
We have
$\Pi_x(f_\beta)/\Pi_x(f_\alpha) > \|\beta\|_0/\|\alpha\|_0$ and
$\Pi_y(f_\beta)/\Pi_y(f_\alpha) < \|\beta\|_0/\|\alpha\|_0$.
Let $\gamma$ and $\delta$ be the strict boundary classes associated
to the ideal points $y$ and $x$ respectively.
Then, using the proof of \cite[Lemma 1.4.1]{cgls},
we see that the number $\Pi_x(f_\beta)/\Pi_x(f_\alpha)$
(resp., $\Pi_y(f_\beta)/\Pi_y(f_\alpha)$)is equal to
$|s_{\tilde \gamma}|$ (resp., $|s_{\tilde \delta}|$).

Next suppose $\Pi_x(f_\alpha)=0$. This happens only in case
$\Pi_x(f_\alpha)/\|\alpha\|_0 < \Pi_x(f_\beta)/\|\beta\|_0$.
Then $\alpha$ is a strict boundary slope and satisfies the desired
condition
$|s_{\tilde\alpha}|=\infty > \|\beta\|_0/\|\alpha\|_0$.

In case $\alpha =\mu$ and $\beta=a\mu +b\lambda$,
by changing coordinates we have $s_\gamma = b(r_\gamma - r_\beta)$.
Hence we have the conclusions of Theorem \ref{main}(1).

Suppose (2) of Proposition \ref{prop} holds.
Note that there is an ideal point $x$
such that $\Pi_x(f_\alpha)>0$, otherwise  $\|\alpha\|_0=0$ and
this is a contradiction to the definition of the Culler-Shalen norm.
Hence for the strict boundary slope $\gamma$ associated to
the ideal point $x$,
we have $|s_{\tilde \gamma}|={\|\beta\|_0}/{\|\alpha\|_0}$.
There are at least two distinct strict boundary classes
associated to $\mathcal I$ by Lemma \ref{two}.
Hence there are exactly two distinct boundary classes,
say $\gamma$ and $\delta$,
such that
$-s_{\tilde \gamma}=\|\beta\|_0/\|\alpha\|_0$ and
$s_{\tilde \delta}=\|\beta\|_0/\|\alpha\|_0$.
(Here we assumed without loss of generality that
$r_\delta>0$.)
In case $\alpha =\mu$ and $\beta=a\mu +b\lambda$,
by changing coordinates we have the conclusion.

The proof for the canonical norm when $M$ is hyperbolic is
exactly the same.
\qed
\vspace{3mm}

Next we prove Corollary~\ref{app}. We will prove the three parts 
separately. In each case we calculate the ratio
$t=\|\beta\|_0/\|\alpha\|_0=\|\beta\|_M/\|\alpha\|_M$
and apply Theorem \ref{main} and Corollary \ref{general}.

First we remark that Theorem~\ref{main}(2) does not occur
when there are two distinct cyclic classes.
Indeed, Dunfield proved the following result.

\begin{lem}[Lemma 4.4 and 4.5\cite{d1}]\label{dunfield}
Suppose $M$ is hyperbolic and $H^1(M,\mathbb Z_2)=\mathbb Z_2$.
Let $\alpha$ and $\beta$ be cyclic classes.
Then $f_\beta/f_\alpha$ cannot be constant on $X_0$.
\end{lem}

He then proved that $|s_{\tilde \gamma}|<1$. Our
Corollary~\ref{app}(1) asserts additionally the opposite inequality
$|s_{\tilde \delta}|>1$.

\vspace{3mm}

\noindent
{\it Proof of Corollary $\ref{app} (1)$}.
If $\beta$ has minimal norm, then $t \leq 1$ and 
we have a proof of the first assertion.
If $\alpha$ and $\beta$ are both cyclic, then $t = 1$.
Due to Proposition~1.1.3 of~\cite{cgls},
the function $f_\beta/f_\alpha$ cannot have poles except at ideal points.
Hence if (2) of Theorem~\ref{main} occurs then
the function $f_\beta/f_\alpha$ also
has no poles at the ideal points and is, therefore, constant.
However this contradicts Lemma \ref{dunfield}.
Thus (1) of Theorem \ref{main} holds.
Hence, we can find two distinct strict boundary classes $\gamma$ and $\delta$
with $|s_{\tilde \gamma}|<1$ and $|s_{\tilde \delta}|> 1$.
\qed
\vspace{3mm}

For the finite slope case,
we quote a lemma of \cite{bz2}.

\begin{lem}[Lemma 9.1\cite{bz2}]\label{s4}
If $M$ is hyperbolic,
$4\leq 2|H_1(M,\mathbb Z_2)|\le s_M$ holds.
\end{lem}

\noindent
{\it Proof of Corollary $\ref{app} (2)$}.
If $\beta$ is a strict boundary class,
then $\beta$ satisfies the 
conclusion.
Hence we assume that $\beta$ is not a strict boundary class.
First note that for a sublattice $\tilde L$ of index $q$,
$q\alpha$ is in $\tilde L$ for any element $\alpha \in L$ and
if $\beta$ realizes the minimal norm on a $\tilde L$,
then $\|\beta\|_M \leq \|q\alpha\|_M=q\|\alpha\|_M$.
Hence we have $t\leq q$.
By Lemma \ref{s4} we have $s_M\geq 4$.
Then by using this fact and Theorem~\ref{finite}:
if $\beta$ is C-type, $t \leq 1$;
if $\beta$ is D-type, $t\leq 2$;
if $\beta$ is T-type, from $\|\beta\|_M\leq s_M+4$ we have $t\leq 2$;
if $\beta$ is I-type, from $\|\beta\|_M\leq s_M+8$  we have $t\leq 3$;
if $\beta$ is O-type and $H_1(M,\mathbb Z)$ has no non-trivial
even torsion,
from $\|\beta\|_M\leq s_M+6$ we have $t\leq 5/2$;
finally, if $\beta$ is O-type and $H_1(M,\mathbb Z)$ has non-trivial
even torsion, then $t \leq q \leq 3$.
Hence we have $t={\|\beta\|_M}/{\|\alpha\|_M}\leq 3$.
\qed
\vspace{3mm}

\noindent
{\it Proof of Corollary $\ref{app} (3)$}.
We assume that $\beta$ is not a strict boundary class.
By Theorem \ref{bb}, we have $\|\beta\|_0=s_0+2A$.
Hence we have $t\leq 1+2A/s_0$.
\qed
\vspace{3mm}

\noindent
{\it Proof of Corollary $\ref{diam}$}.
We take $\alpha$ to be $\mu$
and $\beta$ to be a strict boundary slope associated to an ideal point,
say $x$.
Since $\infty \notin \mathcal B$, we have $\Pi_y(f_\mu)>0$
for any ideal point $y \in \mathcal I$.
Since $\Pi_x(f_\beta)=0$, we have
$\Pi_x(f_\beta)/\|\beta\|_0 < \Pi_x(f_\mu)/\|\mu\|_0$,
i.e., case (1) in Theorem~\ref{main} always holds.
Then we can find a strict boundary class $\delta$ with
$|r_\delta -r_\beta|> {\|\beta\|_0}/{b\|\mu\|_0}$.
Since {\rm diam}\,$\mathcal B \geq |r_\delta -r_\beta|$,
we have the conclusion.
\qed

\section{Examples \label{secex}}
Corollary~\ref{app} shows that a cyclic, finite, or
Seifert slope lies near a strict boundary slope $r_\gamma$.
We verify this conclusion for the twist knots
and the $(-2,3,n)$ pretzel knots by taking $\alpha$  to be the meridian
and $\beta$ to be one of these exceptional classes. We
will also verify the second assertion of Theorem~\ref{main}(1) by
evaluating $t = \| \beta \|_0 / \| \alpha \|_0$ and observing that
$|r_\gamma - r_\beta| < t$. (The exceptional slopes are all integral,
so that $b=1$.)

For each of these knots, there is only one norm curve
$X_0$ in the character variety. Moreover, with the exception of the
figure eight knot (which is a kind of twist knot) the ideal points of
$\tX_0$ are associated to three different strict boundary slopes. This
means that (1) of Proposition~\ref{prop} holds, since (2)
would imply that there are exactly two distinct strict boundary slopes
associated to the norm curve (see Theorem~\ref{main}). For each knot
we determine the strict boundary slopes
associated to the ideal points
$x$ and $y$ of Proposition~\ref{prop}(1).

In addition, we calculate the diameter of $\mathcal B$ for each knot
and compare it with the best estimate obtained from Corollary~\ref{diam}.

\subsection{The $(-2,3,n)$ pretzel knots}

We will assume $n$ is odd and $n \neq 1,3,5$ so that the $(-2,3,n)$
pretzel knot is hyperbolic. The Culler-Shalen seminorms of this knot
are worked out explicitly in~\cite{m1} where it is shown that there
is only one norm curve $X_0$ in the character variety. We first examine
the $(-2,3,7)$ and $(-2,3,9)$ knots which have cyclic and finite slopes
before turning to the remaining pretzel knots which have Seifert slopes.

\subsubsection{The $(-2,3,7)$ pretzel knot}

The finite surgeries of the $(-2,3,7)$ pretzel knot
are classified in~\cite[Example 10.1]{bz1}. There are
cyclic surgeries on the meridian, and at slopes 18 and 19,
as well as an I-type finite surgery at slope 17. The boundary slopes
are given in~\cite{ho} as $0$, $16$, $37/2$, and $20$.
The longitude $0$ is the slope of a fibre in a fibration\cite{gabai} while
the remaining boundary slopes correspond to ideal points of $\tX_0$ and
are therefore strict. Indeed, the calculation of~\cite{m1} shows
that there are ideal points $u$ and
$v$ of $\tX_0$ associated to the slopes $16$ and $20$
respectively with $\Pi_u(f_{\mu}) = \Pi_v(f_{\mu}) = 2$.

For the boundary slope $37/2$, the result
of~\cite{m1} does not allow us to determine the number of associated
ideal points. To calculate the norm of a class $\gamma$, we need only
the total $\Pi_x(f_\gamma)$ summed over all ideal points. Thus,
although we know from~\cite{m1} that $\Pi_x(f_\mu)$ summed over the ideal
points associated with $37/2$ is eight, we cannot
distinguish between two possibilities:
either there is one associated ideal point $w$ with $\Pi_w(f_{\mu}) = 8$
or else there are two ideal points $w_1$, $w_2$, each having
$\Pi_{w_i}(f_{\mu}) = 4$ (by~\cite[Lemma 6.2(1)]{bz1}, $ 4 \mid
\Pi(f_\mu)$ for any $37/2$ ideal point). Since it will not
affect our conclusions below, for the sake of brevity we will assume
that there is only one ideal point $w$.

Given $\Pi_x(f_{\mu})$, the order of pole of any other $f_{\gamma}$ is
determined by the formula (see~\cite[Lemma 6.2(1)]{bz1})
\begin{equation}
  \Pi_x(f_{\gamma}) = \frac{\Delta(r_\gamma, r_\beta)}{\Delta(r_\mu,
r_\beta)}\Pi_x(f_{\mu}), \label{eq621}
\end{equation}
where $\beta$ is the boundary class
associated to $x$ and
$\Delta(a/b,c/d) = |ad-bc|$ is the minimal geometric intersection of the
two slopes. For example, Table~\ref{tbl1} gives the degree of pole of
various functions at the three ideal points.
\begin{table}[htbp]
\begin{tabular}{c|c|cccc}
Ideal& Associated & $\Pi(f_{\mu})$ & $\Pi(f_{17})$ &
$\Pi(f_{18})$ &
$\Pi(f_{19})$ \\
point & boundary slope \\ \hline
$u$ & $16$ &  $2$ & $2$ & $4$ & $6$ \\
$v$ & $20$ & $2$ & $6$ & $4$ & $2$ \\
$w$ & $37/2$ & $8$ & $12$ & $4$ & $4$
\end{tabular}
\caption{Order of pole at ideal points of $\tX_0$ for the $(-2,3,7)$
pretzel knot \label{tbl1}}
\end{table}

Let
$\alpha = \mu$ (then~\cite{m1}, $\|\alpha\|_0 = s_0 = 12$) and let $\beta$
be one of the other cyclic or finite classes. We will determine the
boundary slope associated to the ideal points $x$ and $y$ of
Proposition~\ref{prop}(1).

If $\beta = 19$, $\|\beta\|_0 = 12$ so that $t = \| \beta \|_0 / \|
\alpha \|_0 = 1$. Here,
$\Pi_u(f_{\alpha})/\|\alpha\|_0 <
\Pi_u(f_{\beta})/\|\beta\|_0$ while $\Pi_w(f_{\alpha})/\|\alpha\|_0 >
\Pi_w(f_{\beta})/\|\beta\|_0$. Thus, in the Proposition, $x = u$ and
$y = w$. Moreover, the boundary slopes of Corollary~\ref{app} are
the associated boundary slopes $r_\gamma = 37/2$ and $r_\delta = 16$.
These also verify the second assertion of Theorem~\ref{main}(1) since
$| r_\gamma - r_\beta | = |37/2 - 19| = 1/2 < t = 1 $ and
$| r_\delta - r_\beta | = |17 - 19| = 2 > t = 1 $.

For $\beta = 18$, again, $\|\beta\|_0 = 12$  and $t = 1$. Here,
$\Pi_u(f_{\alpha})/\|\alpha\|_0 <
\Pi_u(f_{\beta})/\|\beta\|_0$
and $\Pi_v(f_{\alpha})/\|\alpha\|_0 <  \Pi_v(f_{\beta})/\|\beta\|_0$,
while $\Pi_w(f_{\alpha})/\|\alpha\|_0 >
\Pi_w(f_{\beta})/\|\beta\|_0$. Therefore, in Corollary~\ref{app}, we
again have $r_\gamma=37/2$ while $16$ and $20$ are both valid choices for
$r_\delta$. That is, $|r_\gamma - r_\beta| = 1/2 < t =1$ and
$|r_\delta - r_\beta| = 2 > t = 1$.

Finally, $\beta = 17$ has norm $\|\beta\|_0 = 20$~\cite{m1} so that
$t = 20/12 = 5/3$. Here,
$\Pi_v(f_{\alpha})/\|\alpha\|_0 <
\Pi_v(f_{\beta})/\|\beta\|_0$
while $\Pi_u(f_{\alpha})/\|\alpha\|_0 > \Pi_u(f_{\beta})/\|\beta\|_0$
and $\Pi_w(f_{\alpha})/\|\alpha\|_0 >
\Pi_w(f_{\beta})/\|\beta\|_0$. Hence, $37/2$ and $16$ are strict
boundary slopes $r_\gamma$ near the finite slope $r_\beta = 17$. Note
that $| r_\gamma - r_\beta | < 5/3 = t$ in both cases. In other words,
for finite slopes, $t/b$ will often give us a better estimate than the
bound of $3/b$ stated in  Corollary~\ref{app}.

The diameter of the set of strict boundary slopes is $20 - 16 = 4$.
Using Corollary~\ref{diam}, we obtain the lower bound
$\| 16 \|_0 / \| \mu \|_0 = 28/12 = 7/3$.

\subsubsection{The $(-2,3,9)$ pretzel knot}

The finite surgeries of the $(-2,3,9)$ pretzel knot
are classified in~\cite{m1}. There's a cyclic meridional surgery,
an O-type finite surgery at slope 22, and an I-type finite surgery at
slope 23. The boundary slopes may be calculated using the
methods of~\cite{ho,d2} as $0$, $16$, $67/3$ and $24$.
Only $0$ is not strict~\cite{gabai}.

The calculation of~\cite{m1} shows that there are
ideal points $u$ and $v$ associated to the slopes $16$ and $24$.
Again the situation for $67/3$ is left ambiguous. Although there may be
two ideal points each with $\Pi(f_\mu) = 6$, we will simply assume
that there is only one $67/3$ ideal point $w \in \tX_0$ with
$\Pi_w(f_\mu) = 12$. The order of pole for various $f_{\gamma}$ is shown
in Table~\ref{tbl2}.
\begin{table}[htbp]
\begin{tabular}{c|c|ccc}
Ideal & Associated & $\Pi(f_{\mu})$ & $\Pi(f_{22})$ & $\Pi(f_{23})$  \\
point & boundary slope \\ \hline
$u$ & $16$ & $2$ & $12$ & $14$ \\
$v$ & $24$ & $2$ & $4$ & $2$  \\
$w$ & $67/3$ & $12$ & $4$ & $8$
\end{tabular}
\caption{Order of pole at ideal points of $\tX_0$ for the
$(-2,3,9)$ pretzel knot \label{tbl2}}
\end{table}

Let $\alpha = \mu$, $\beta_1 = 22$ and $\beta_2
  = 23$.
Then~\cite{m1}, $\| \alpha \|_0 = 16$, $\| \beta_1 \|_0 = 20$ and
$\| \beta_2 \|_0 = 24$ so that $t_1 = \| \beta_1 \|_0 / \| \alpha \|_0 =
5/4$ and $t_2 = 3/2$.

For $\beta_1$,
$\Pi_u(f_{\alpha})/\|\alpha\|_0 <
\Pi_u(f_{\beta_1})/\|\beta_1\|_0$ and
$\Pi_v(f_{\alpha})/\|\alpha\|_0  <
\Pi_v(f_{\beta_1})/\|\beta_1\|_0$ while
$\Pi_w(f_{\alpha})/\|\alpha\|_0  >
\Pi_w(f_{\beta_1})/\|\beta_1\|_0$. In other words, we can choose
$y = w$ in Proposition~\ref{prop} while $x = u$ and $x = v$ are
both valid choices. In Corollary~\ref{app}, we have $r_\gamma = 67/3$.
Note that $|r_\gamma - r_{\beta_1}| = |67/3 - 22| = 1/3 < 5/4 = t_1$.

For $\beta_2$, $x=u$ while $y = v$ and $y=w$ are both correct
in Proposition~\ref{prop}. Consequently, $r_\gamma
= 67/3$ and $r_\gamma = 24$ both satisfy Corollary~\ref{app}. Again, these
in fact satisfy the stronger inequality $|r_\gamma - r_{\beta_2}| < 3/2
= t_2$.

The diameter of the set of strict boundary slopes is $24 - 16 = 8$.
Using Corollary~\ref{diam}, we obtain the lower bound
$\| 16 \|_0 / \| \mu \|_0 = 92/16 = 23/4$.

\subsubsection{Pretzel
 knots with Seifert slopes}

For $n$ odd and $n \neq 1,3,5,7,9$, the $(-2,3,n)$ pretzel knot admits
Seifert surgeries at slopes $\beta_1 = 2n+4$ and $\beta_2 = 2n+5$
(see~\cite{bh}). The number of ideal points may be quite large for these
knots, so we will work with a slight reformulation of
Propostion~\ref{prop}.  Using Equation~\ref{eq621},
\begin{equation}
\Pi_x(f_\alpha)/\| \alpha \|_0 > \Pi_x(f_\beta) / \| \beta\|_0
\mbox{ if and only if } \Delta(r_\alpha, r_\delta) / \| \alpha \|_0 >
\Delta(r_\beta, r_\delta) / \| \beta \|_0 \label{eqprop71}
\end{equation}
where $\delta$ is the strict boundary class associated to the ideal point
$x$. The boundary slopes and norm of the $(-2,3,n)$ pretzel knot differ
depending on the sign of $n$, so we consider the two cases separately.

If $n \geq 11$, the boundary slopes are~\cite{ho,d2} $0$, $16$, $2n+6$,
and $\frac{n^2-n-5}{\frac{n-3}{2}}$, but $0$ is not strict~\cite{gabai}.
The norm depends on whether or not $3 \mid n$: $\| \alpha \|_0  
= \| \mu \|_0 = s_0 = 3(n-3)$ (respectively, $3n-11$), 
$\| \beta_1 \|_0 = 6(n-5)$ (resp.,
$6n-34$), and $\| \beta_2 \|_0 = 7n-37$ (resp., $7n-39$) when $3
  {\not |~} n$ (resp., $3 \mid n$). So, $t_1 = 2 (n-5)/ (n-3)$
(respectively,
$\frac{6n-34}{3n-11}$) and $t_2 = \frac{7n-37}{3(n-3)}$ (resp.,
$\frac{7n-39}{3n-11}$). Thus, $\Delta(r_\mu, r_\delta) / \| \mu \|_0 >
\Delta(r_{\beta_1}, r_\delta)/ \| \beta_1 \|_0$ only for the boundary
slope $\delta = \frac{n^2-n-5}{\frac{n-3}{2}}$.

Indeed, the number of non-abelian characters for this Seifert slope is
$A_1 = \frac32 (n-7)$ (respectively, $\frac12 (3n-23)$) when $3 {\not |~}
n$ (resp., $3 \mid n$), and, for $r_\gamma =
\frac{n^2-n-5}{\frac{n-3}{2}}$, $| r_\gamma - r_{\beta_1}| = 2/(n-3) <
t_1 = 1 + 2A_1/s_0$ in agreement with Corollary~\ref{app}.
For the other Seifert slope, $\beta_2 = 2n+5$, the inequality
of Proposition~\ref{prop}, $\Delta(r_\mu, r_\delta) / \| \mu \|_0 > 
\Delta(r_{\beta_2}, r_\delta)/ \| \beta_2 \|_0$, holds for both
$\delta = 2n+6$ and $\frac{n^2-n-5}{\frac{n-3}{2}}$. Indeed, the number
of non-abelian characters is $A_2 = 2(n-7)$ and $|r_\gamma -
r_{\beta_2}| \leq 1 < t_2 = 1 + 2A_2/s_0$ for both $r_\gamma = 2n+6$ and
$r_\gamma = \frac{n^2-n-5}{\frac{n-3}{2}}$.

When $n \geq 11$, the diameter of the set of strict boundary slopes 
of the $(-2,3,n)$ pretzel knot is $2n+6- 16 = 2(n-5)$. 
Using Corollary~\ref{diam}, we obtain the lower bound
$\| 16 \|_0 / \| \mu \|_0 = (6n^2-56n+126)/(3n-9)$ (respectively,
$(6n^2 - 60n + 146)/(3n-11)$) when $3 {\not |~} n$ (resp., $3 \mid n$).
Thus, the difference, (diam\,$\mathcal B - \| 16 \|_0 / \| \| \mu \|_0$),
is $13/6$ when $n = 11$ and increases towards $8/3$ as $n$ tends to 
infinity.

If $n < 0$, the boundary slopes are ~\cite{ho,d2} $0$, $10$, $2n+6$, and
$\frac{2(n+1)^2}{n}$. The longitude $0$ is not strict~\cite{gabai} unless
$n = -1$ or $n = -3$. Again, the norm depends on whether or not $3 \mid n$
(see~\cite{m1}): $\| \mu \|_0 = s_0 = 3(1-n)$ (respectively, $1-3n$),
$\| \beta_1 \|_0 = 6(3-n)$ (resp., $2(7-3n)$), and
$\| \beta_2 \|_0 = 15-7n$ (resp., $13-7n$) when $3 {\not |~} n$ (resp.,
$3 \mid n$). So, $t_1 = 2 (3-n)/(1-n)$ (respectively, $2(7-3n)/(1-3n)$)
and $t_2 = \frac{15-7n}{3(1-n)}$ (resp., $\frac{13-7n}{1-3n}$).
Thus, both $\delta = 2n+6$ and $\delta =
\frac{2(n+1)^2}{n}$ will satisfy the Proposition~\ref{prop} inequality
(see Equation~\ref{eqprop71})
$\Delta(r_\mu, r_\delta) / \| \mu \|_0 >
\Delta(r_\beta, r_\delta)/ \| \beta \|_0$ for the Seifert slopes $\beta_1
= 2n+4$ and $\beta_2 = 2n+5$. Indeed, for $\beta_1$, the number of
non-abelian characters is $A_1 = \frac32 (5-n)$ (respectively, $\frac12
(13-3n)$) when $3 {\not |~} n$ (resp., $3 \mid n$) and
$|r_\gamma -r_{\beta_1} | < t_1
= 1 + 2A_1/s_0$ for both $r_\gamma = 2n+6$ and
$r_\gamma = \frac{2(n+1)^2}{n}$. For $\beta_2$, we have $A_2 =2(3-n)$
and, again, both choices of $r_\gamma$ verify Corollary~\ref{app}:
$|r_\gamma - r_{\beta_2}| < t_2 = 1 + 2A_2/s_0$.

When $n < 0$, the diameter of the set of strict boundary slopes 
of the $(-2,3,n)$ pretzel knot is $10 - 2(n+1)^2/n = 6 -2n-2/n$. 
Using Corollary~\ref{diam}, we obtain the lower bound
$\| 10 \|_0 / \| \mu \|_0 = (6n^2-18n+8)/(3-3n)$ (respectively,
$(6n^2 - 14n)/(1-3n)$ when $3 {\not |~} n$ (resp., $3 \mid n$).
Thus, the difference, (diam\,$\mathcal B - \| 10 \|_0 / \| \| \mu \|_0$),
is $14/3$ when $n = -1$ and decreases towards $2$ as $n$ tends to 
negative infinity.

\subsection{The twist knot $K_n$}

\begin{figure}
\begin{center}
\includegraphics{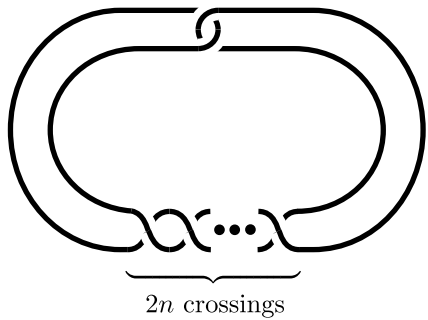}
\caption{\label{Kn}%
The twist knot $K_n$.}
\end{center}
\end{figure}

Figure~\ref{Kn}
shows the twist knot $K_n$.
We will assume $n \neq 0,1$ so that the $K_n$
twist knot is hyperbolic.
These knots have Seifert slopes at $-1$, $-2$, and $-3$.
Burde~\cite{b} showed that the character
variety has only one norm curve $X_0$ and the associated
Culler-Shalen seminorm is determined in~\cite{bmz}.
Ohtsuki~\cite{o} has enumerated the ideal points $x$ of these knots and
demonstrated that $\Pi_x(f_{\mu}) = 2$ at each ideal point.
Since the norm and boundary slopes depend on the sign of $n$, we consider
two cases.

If $n \geq 2$, the boundary slopes are~\cite{ht} $0$, $-4$, and
$-(4n+2)$ and these are all strict~\cite{bmz,o}. The norms are
$\| \mu \|_0 = s_0 = 4n-2$, $\| -1 \|_0 = 2(8n-3)$, $\|-2\|_0 =
8(2n-1)$, and $\|-3\|_0 = 2(8n-5)$. For each of the Seifert slopes
$\beta$, the inequality of Proposition~\ref{prop} 
(see Equation~\ref{eqprop71}),
$\Delta(r_\mu, r_\delta)
/\|\mu\|_0 >
\Delta(r_\beta, r_\delta)/\|\beta\|_0$, obtains when $\delta$ is either of
the boundary slopes $0$ or $-4$. Indeed, the number of non-abelian
characters is $A = 6n-2$ (respectively, $6n-3$, $6n-4$) for the Seifert
slope
$-1$ (resp., $-2$, $-3$) so that $|r_\gamma - r_\beta| \leq t = 1 +
2A/s_0$ whenever $r_\beta$ is one of the three Seifert slopes and
$r_\gamma$ is one of the boundary slopes $0$ or $-4$, in agreement with
Corollary~\ref{app}.

The diameter of $\mathcal B$ is $4n+2$ when $n \geq 2$. 
Using Corollary~\ref{diam}, we obtain the lower bound
$\| -(4n+2) \|_0 / \| \mu \|_0 = 16n(n-1)/(4n-2) = 8n(n-1)/(2n-1)$.
The difference between the diameter and this bound is $14/3$ when
$n = 2$ and decreases towards $4$ as $n$ tends to infinity.

If $n \leq -1$, the boundary slopes are~\cite{ht} $0$, $-4$ and $-4n$ and
these are strict as long as $n \leq -2$. For the figure eight knot,
$K_{-1}$, $0$ is not a strict boundary slope (but $\pm 4$ are).
The minimal norm is~\cite{bmz} $\| \mu \|_0 = -4n$ and the Seifert slopes
$-1$, $-2$, and $-3$ all have norm $-16n$. Again, the boundary slopes
$0$ and $-4$
satisfy the inequality of Propostion~\ref{prop},
$\Delta(r_\mu, r_\delta) /
\|\mu\|_0 >
\Delta(r_\beta, r_\delta)/\|\beta\|_0$, for each of the Seifert slopes
$\beta$. Indeed, the number of non-abelian characters is $A = -6n$ for
each of the three Seifert slopes so that $|r_\gamma - r_\beta| \leq 4 = t
= 1+2A/s_0$, in accord with Corollary~\ref{app}, whenever $r_\beta$ is
Seifert and $r_\gamma$ is one of the boundary slopes $0$ or $-4$.

The diameter of $\mathcal B$ is $4-4n$ when $n \leq -1$. 
Using Corollary~\ref{diam}, we obtain the lower bound
$\| -4n \|_0 / \| \mu \|_0 = 16n^2/(-4n) = -4n$.

The figure eight knot, $K_{-1}$, is of special interest as it provides an
example of Theorem~\ref{main}(2) and Proposition~\ref{prop}(2).
For this knot, the norm
curve $\tX_0$ has only two associated strict boundary slopes $4$ and $-4$.
Let
$r_\alpha = a/b$ and $r_\beta = c/d$. Then $\alpha$ and $\beta$ will
satisfy  part 2 of Theorem \ref{main} and Proposition \ref{prop}
provided $16bd = ac$. For
example, let $r_\alpha = 1/0$ (so that $\alpha = \mu$) and $r_\beta = 0/1$
($\beta = \lambda$).
Then, $\| \beta \|_0 / \| \alpha \|_0 =
 16/4 = 4$, so that $ r_{\beta} - r_{\gamma} = \|\beta\|_0/\|\alpha\|_0$
for $r_{\gamma} = -4$, and $r_{\delta} - r_{\beta} = 
\|\beta\|_0/\|\alpha\|_0$ for $r_{\delta} = 4$
 (compare Theorem~\ref{main}(2)).
For ideal points $u$ associated to the slope $4$, we have
$\Delta(r_\alpha, 4)/ \| \alpha \|_0 = 1/4$ and $ \Delta(r_\beta,4)/ \|
\beta \|_0 = 4/16 = 1/4$. Therefore, (compare Equation~\ref{eqprop71})
$\Pi_u(f_\alpha)/\|\alpha\|_0 = \Pi_u(f_\beta) \|\beta\|_0$ in accord
with Proposition~\ref{prop}(2). Similarly, at any ideal point $v$
associated to slope $-4$, $\Pi_v(f_\alpha) / \| \alpha \|_0 =
\Pi_v(f_\beta) / \| \beta \|_0$ since
$\Delta(r_\alpha,-4)/ \| \alpha \|_0 = 1/4$ and
$\Delta(r_\beta,-4) / \| \beta \|_0 = 4/16 = 1/4$.
\vspace{2cm}

\noindent
{\bf Acknowledgments.}

This paper is based on research begun during a
visit of the second author to Tohoku University. He
wishes to thank Prof.\ Hajime Urakawa and the Graduate School of
Information Sciences for the invitation and hospitality
during his stay.
The second author also wishes to thank Steve Boyer for
helpful conversations.
The third author would like to express his gratitude
to Prof. Daryl Cooper for precious comments.


\begin{thebibliography}{CGLS}

\bibitem[BB]{bb}
     L. Ben Abdelghani and S. Boyer,
     A calculation of the Culler-Shalen seminorms associated to
     small Seifert Dehn fillings,
     Proc. London Math. Soc.
     \textbf{83} (2001), 235--256.

\bibitem[BH]{bh}
    S. Bleiler and C. Hodgson,
     Spherical space forms and Dehn fillings,
     Topology
     \textbf{35} (1996), 809--833.

\bibitem[BMZ]{bmz}
     S. Boyer, T.W. Mattman, and X. Zhang,
     The fundamental polygons of twist knots and the $(-2,3,7)$
     pretzel knot,
     KNOTS '96 (Tokyo),  World Sci. Publishing,
     (1997), 159--172.

\bibitem[BZ1]{bz1}
     S. Boyer and X. Zhang,
     Finite Dehn surgery on knots,
     J. Amer. Math. Soc.
     \textbf{9} (1996), 1005--1050.

\bibitem[BZ2]{bz2}
     S. Boyer and X. Zhang,
     A proof of the finite filling conjecture,
     J. Diff. Geom.
     \textbf{59} (2001), 87--176.

\bibitem[B]{b}
      G. Burde,
      ${\rm SU}(2)$-representation spaces for
      two-bridge knot groups,
      Math. Ann.
      \textbf{288} (1990), 103--119.

\bibitem[CGLS]{cgls}
      M. Culler, C.McA. Gordon, J. Luecke and P.B. Shalen,
      Dehn surgery on knots,
      Ann. of Math.
      \textbf{125} (1987), 237--300.

\bibitem[CS1]{cs}
      M. Culler and P.B. Shalen,
      Varieties of group representation and splittings of $3$-manifolds,
      Ann. of Math.
      \textbf{117}(1983), 109--146.

\bibitem[CS2]{cs2}
      M. Culler and P.B. Shalen,
      Bounded, separating surfaces in knot manifolds,
      Invent. Math.
      \textbf{75} (1984), 537--545.

\bibitem[CS3]{cs3}
      M. Culler and P.B. Shalen,
      Boundary slopes of knots,
      Comment. Math. Helv.
      \textbf{74} (1999), 530--547.

\bibitem[D1]{d1}
      N. Dunfield,
      Cyclic surgery, degrees of maps of character curves, and volume
      rigidity for hyperbolic manifolds,
      Invent. Math.
      \textbf{136} (1999), 623--657 math.GT/9802022.

\bibitem[D2]{d2}
      N. Dunfield,
      A table of boundary slopes of Montesinos knots,
      Topology
      \textbf{40} (2001), 309--315 math.GT/9901120.

\bibitem[G]{gabai}
       D. Gabai,
       Detecting fibred links in $S^3$,
       Comment.\ Math.\ Helvetici
       \textbf{61} (1986), 519--555.

\bibitem[HO]{ho}
      A. Hatcher and U. Oertel,
      Boundary slopes for Montesinos knots,
      Topology
      \textbf{28} (1989), 453--480.

\bibitem[HT]{ht}
    A. Hatcher and W. Thurston,
    Incompressible surfaces in 2-bridge knot complements,
    Invent. Math.
    \textbf{79} (1985), 225--246.

\bibitem[K]{kirby}
      R. Kirby,
      Problems in low-dimensional topology,
      Geometric Topology, Volume 2, editor W. Kazez,
      AMS/IP Studies in Advanced Mathematics, 1996.

\bibitem[Ma]{m1}
      T.W. Mattman,
      The Culler-Shalen seminorms of the $(-2,3,n)$ pretzel knot,
      J. Knot Theory and Ram. 
      \textbf{11} (2002), 1251--1289 math.GT/9911085.

\bibitem[Mi]{mi}
      J. Milnor,
      Groups which act on $S^n$ without fixed points,
      Amer. J. Math.
      \textbf{79} (1957), 623--631.

\bibitem[O]{o}
      T. Ohtsuki,
      Ideal points and incompressible surfaces in two-bridge knot
      complements,
      J. Math. Soc. Japan
      \textbf{46} (1994) 51--87.

\end{thebibliography}
\end{document}